\newif\ifspringer
\IfFileExists{svmult-5.10.cls}{\springertrue}{\springerfalse}
\ifspringer
\documentclass[graybox]{svmult}

\usepackage{type1cm}        % activate if the above 3 fonts are
\usepackage{makeidx}         % allows index generation

\usepackage[bottom]{footmisc}% places footnotes at page bottom

\usepackage{newtxtext}       % 
\usepackage[varvw]{newtxmath}       % selects Times Roman as basic font

\makeindex             % used for the subject index
\else

\documentclass[twoside,a4paper]{amsart}
\IfFileExists{eulervm.sty}{%
\usepackage{eulervm}}
{
}

\fi 

\usepackage{bbm}
\providecommand{\sperp}{\circ}
\usepackage{mathtools}
\DeclareFontEncoding{LS1}{}{}
\DeclareFontSubstitution{LS1}{stix2}{m}{n}
\DeclareSymbolFont{symbols2}{LS1}{stix2frak}{m}{n}
\DeclareMathSymbol{\obot}{\mathbin}{symbols2}{"A9}
\providecommand{\operp}{\obot}
\providecommand{\Lin}{\mathrm{Lin}}

\PassOptionsToPackage{pdfauthor={afar Aljasem and Vladimir V. Kisil},%
    pdftitle={Operator Projective Line and Its Transformations},%
    pdfsubject={mathematics},%
    backref=page,%
  pdfkeywords={symmetries}
}{hyperref}

\providecommand{\Space}[3][]{\ifx#2R\ifx#1e \mathbb{C}^{#3} \else
\ifx#1p \mathbb{D}^{#3} \else
\ifx#1h \mathbb{O}^{#3} \else
\ifx#1\sigma \mathbb{A}^{#3} \else
\ensuremath{\mathbb{#2}^{#3}_{#1}{}} \fi \fi \fi \fi \else
\ensuremath{\mathbb{#2}^{#3}_{#1}{}} \fi}
  \providecommand{\FSpace}[3][]{\ensuremath{\ifx#2l \ell_{#3}^{#1}{}\else
  \mathsf{#2}_{#3}^{#1}{}\fi}}

\providecommand{\SL}[1][2]{\ensuremath{\mathrm{SL}_{#1}(\Space{R}{})}}
\providecommand{\Sp}[1][n]{\ensuremath{\mathrm{Sp}(#1)}}

\providecommand{\oper}[1]{\mathcal{#1}}
\providecommand{\such}{\,\mid\,}
\providecommand{\rme}{\mathrm{e}}
\providecommand{\norm}[2][\relax]{\left\|#2\right\|\ifx#1\relax\else_{#1}\fi}
\providecommand{\modulus}[2][\relax]{\left| #2 \right|\ifx#1\relax\else_{#1}\fi}

\providecommand{\scalar}[3][\relax]{\left\langle #2,#3 
        \right\rangle\ifx#1\relax\else_{#1}\fi}
        \providecommand{\MR}[1]{MR\href{http://www.ams.org/mathscinet-getitem?mr=#1}{#1}}
  \providecommand{\Zbl}[1]{Zbl\href{http://www.emis.de:80/cgi-bin/zmen/ZMATH/en/zmathf.html?first=1&maxdocs=3&type=html&an=#1&format=complete}{#1}}

  \usepackage[breaklinks=true,%
  colorlinks=true,%
  bookmarks=true,%
  backref=page,%
  pagebackref=true
]{hyperref}
  \AtBeginDocument{\hypersetup{breaklinks=true,%
      colorlinks=true,%
      bookmarks=true,%
      backref=page,%
      pagebackref=true
}}

\ifspringer
  \providecommand{\cites}[1]{\cite{#1}}        
\else
\usepackage[nobysame]{amsrefs}
  \theoremstyle{theorem}
  \newtheorem{theorem}{Theorem}
  \newtheorem{lemma}[theorem]{Lemma}
    \newtheorem{proposition}[theorem]{Proposition}

  \theoremstyle{definition}
    \newtheorem{definition}[theorem]{Definition}
    \newtheorem{example}[theorem]{Example}

    \theoremstyle{remark}
    \newtheorem{remark}[theorem]{Remark}
\fi

\providecommand{\Nodin}{{\#}}
\providecommand{\Nnull}{{\bullet}}

\begin{document}
\newcommand{\myabstract}{We introduce a concept of the operator (non-commutative) projective line \(P\FSpace{H}{}\) defined by a Hilbert space \(\FSpace{H}{}\) and a symplectic structure on it. Points of \(P\FSpace{H}{}\) are Lagrangian subspaces of \(\FSpace{H}{}\). If a particular Lagrangian subspace is fixed then we can define \(\SL\)-action on \(P\FSpace{H}{}\). This gives a consistent framework for linear fractional transformations of operators. Some connections with spectral theory are outline as well.
}

\ifspringer
\title*{Operator Projective Line and \\Its Transformations}
\author{Jafar Aljasem and\\ Vladimir V. Kisil\orcidID{0000-0002-6593-6147}}
\institute{Jafar Aljasem \at University of Leeds, Leeds LS2\,9JT, UK
\email{ml17jhaa@leeds.ac.uk}
\and Vladimir V. Kisil \at University of Leeds, Leeds LS2\,9JT, UK \email{V.Kisil@leeds.ac.uk}}
\maketitle

\centerline{\emph{For Ilya Spitkovsky on his 70th birthday}\medskip}

\abstract{\myabstract}

\else

\title{Operator Projective Line and Its Transformations}

\author[Jafar Aljasem and Vladimir V. Kisil]%
{Jafar Aljasem and  \href{http://v-v-kisil.scienceontheweb.net/}{Vladimir V. Kisil}}
\thanks{On  leave from Odessa University.}

\address{%            
%Institute of Mathematics\\
%Economics and Mechanics\\
%Odessa State University\\
%ul. Petra Velikogo, 2\\
%Odessa-57, 270057, UKRAINE
School of Mathematics\\
University of Leeds\\
Leeds LS2\,9JT\\
UK
}

\email{\href{mailto:ml17jhaa@leeds.ac.uk}{ml17jhaa@leeds.ac.uk}}
\email{\href{mailto:V.Kisil@leeds.ac.uk}{V.Kisil@leeds.ac.uk}}

\urladdr{\url{http://v-v-kisil.scienceontheweb.net/}}

\date{\today}

\begin{abstract}
  \myabstract
\end{abstract}
\dedicatory{For Ilya Spitkovsky on his 70th birthday}
\maketitle
\fi 

\section{Introduction}
\label{sec:introduction}
\hfill\parbox{.5\textwidth}{\small \emph{The graph of a linear transformation on a Hilbert space is very much like a line in a plane.} \\\centerline{P.\,R. Halmos, \cite{Halmos69a}}}
\medskip

Geometry of two subspaces in a Hilbert space is an astonishingly rich area of research with numerous connections to other topics and applications, see the excellent survey~\cite{BoettcherSpitkovsky10a}. Notably, B\"ottcher and Spitkovsky did not put
there the closing keystone to the completed construction, instead~\cite{BoettcherSpitkovsky10a} serves as a convenient consolidation point for further advances~\cites{BoetcherSimonSpitkovsky17a,BoetcherSpitkovsky18a,BoetcherSpitkovsky21a,BoetcherSpitkovsky23a,Spitkovsky18a} continuing some earlier works~\cites{Spitkovsky94a,VasSpi81}.

% from active contributors~\cites{}
Our interest in two subspaces' arrangement came from works in functional calculus and associated spectra. It was proposed in~\cite{Kisil95i} to define a joint spectrum for several non-commuting operators through linear fractional (M\"obius) transformations of operators. It turned out that the approach is of interest even for finite dimensional operators (matrices)~\cites{Loring15a,DeBonisLoringSverdlov22a,CerjanLoring24a}. Furthermore, the related concept of \emph{covariant spectrum}%
\index{covariant spectrum} produced by linear fractional transformations is meaningful for a single non-normal operator as well~\cites{Kisil02a,Kisil11c}. Such transformations for real numbers
\begin{equation}
  \label{eq:linear-fractional-naive}
  \begin{pmatrix}
    a&b\\c&d
  \end{pmatrix}: \ t \mapsto \frac{at+b}{ct+d},%\qquad \text{ where } t\in \Space{R}{},
\end{equation}
immediately raise an issue if \(ct+d=0\). A popular quick fix looks like this~\cite[\S~3.1]{Kisil12a}:
% \begin{quote}
%   \ldots we agree to assign the symbol \(\infty\) to all fractions \(\frac{a}{0}\) with \(a\neq0\); some rudimentary arithmetic includes rules \(\frac{b}{\infty}=0\), \(a\cdot \infty=\infty\) , \(b+\infty=\infty\), etc. but we do not define values of \(\frac{0}{0}\), \(0\cdot \infty\) or \(\infty-\infty\)\ldots  
% \end{quote}
\begin{quote}
  \ldots to extend arithmetic of
    real numbers \ldots by the
    following rules involving infinity:
    \begin{displaymath}
      \frac{a}{0}=\infty, \qquad
      \frac{a}{\infty}=0, \qquad
      \frac{\infty}{\infty}=1, \quad
      \text{ where } a\neq 0.
    \end{displaymath}
    As a consequence we also note that \(a\cdot\infty+b=a\cdot\infty\)
    for \(a\neq 0\).  But we do not define (and, actually, will not need) values of \(\frac{0}{0}\), \(0\cdot \infty\) or \(\infty-\infty\)\ldots
\end{quote}
Obviously, one can wonder why such a na\"{\i}ve patch does not fall apart in a rigorous mathematical context. The explanation is:  the above suggestion is implicitly based on a well-defined mathematical object---the \emph{projective line}%
\index{projective line} \(P\Space{R}{}\), cf.~\S\,\ref{sec:cycles} below. The principal observation: points of \(P\Space{R}{}\) do not form an algebraic structure, e.g. a field or algebra, but they merely constitute a homogeneous space for the action~\eqref{eq:linear-fractional-naive} of the group \(\SL\).

A situation becomes even more delicate if we try to define linear fractional transformations of an operator \(\oper{T}\) acting on a Hilbert space \(\FSpace{H}{}\) by the direct formula:
\begin{equation}
  \label{eq:linear-fractional-naive}
  \begin{pmatrix}
    a&b\\c&d
  \end{pmatrix}: \ \oper{T} \mapsto ({c\oper{T}+d\mathbbm{1}})^{-1}{(a\oper{T}+b\mathbbm{1})}.
\end{equation}
The simplest suggestion is to avoid values of \(-\frac{d}{c}\) from the spectrum of \(\oper{T}\). However, appearing technicalities significantly exceed ones already met in the scalar case.

This paper presents an extension of the projective line construction to an operator setting.  Recall that points of  \(P\Space{R}{}\)  are represented by lines through the origin in \(\Space{R}{2}\), which can be in turn characterised as Lagrangian (maximal isotropic) subspaces with respect to the standard symplectic structure on \(\Space{R}{2}\). Thus, we define points of the operator projective line to be Lagrangian subspaces of a Hilbert space \(\FSpace{H}{}\) with a symplectic structure.

Selecting a Lagrangian subspace \(\FSpace{M}{}\subset \FSpace{H}{}\) defines ``the horizontal axis'' with \(\FSpace[\perp]{M}{}\) serving as its vertical counterpart. Thereafter, another Lagrangian subspace \(\FSpace{N}{}\subset \FSpace{H}{}\) represents the graph of a \emph{self-adjoint} operator \(\oper{T}: \FSpace{M}{} \rightarrow \FSpace{M}{}\), but if we put few additional conditions on the mutual disposition of \(\FSpace{M}{}\) and \(\FSpace{N}{}\). We prefer to lift any restrictions on \(\FSpace{N}{}\) (except the Lagrangian property) and extend the concept of an operator in such a way that  \(\FSpace{N}{}\) still can be treated as its graph. Such ``operators'' are called bi-Fredholm for the reasons discussed in \S\,5. Then, the linear transformation defined by a matrix \(
\begin{pmatrix}
  a&b\\c&d
\end{pmatrix}
\in\SL
\) through coordinate axes \(\FSpace[\perp]{M}{}\oplus\FSpace{M}{}\)  transforms \(\FSpace{N}{}\) and provides a consistent interpretation of the linear fractional formula~\eqref{eq:linear-fractional-naive} for \(\oper{T}\). Such transformations are closely connected to spectral properties of \(\oper{T}\) as will be indicated in \S\,\ref{sec:sl-action-lagrangian}.

We finally remark that the standard two projections theory is traditionally considered in complex Hilbert spaces~\cite{BoettcherSpitkovsky10a} because spectral theory and functional calculus are heavily employed to make short and elegant proofs~\cite{Halmos69a}. Since our aim~\cite{Kisil02a} is to revise the both---the spectral theory and functional calculus---we are starting from the scratch and allow our Hilbert spaces to be real.

\section{M\"obius transformations and symplectic structure}
\label{sec:cycles}
 The
group \(\SL\) is formed by  real \(2\times 2\) matrices with the unit determinant.
% \(\SL\) acts on the real line by linear-fractional maps:
% \begin{equation}
%   \label{eq:moebius-map-defn}
%   \begin{pmatrix}
%     a&b\\c&d
%   \end{pmatrix}: x \mapsto \frac{ax+b}{cx+d},
%   \quad \text{ where } 
%   x\in\Space{R}{}\text{ and } 
%     \begin{pmatrix}
%     a&b\\c&d
%   \end{pmatrix}\in\SL.
% \end{equation}
It naturally %The group \(\SL\)
acts on \(\Space{R}{2}\) by matrix multiplication on
column vectors:
\begin{equation}
  \label{eq:left-matrix-mult}
  \oper{L}_g:\ 
  \begin{pmatrix}
    x_1\\x_2
  \end{pmatrix}
  \mapsto
  \begin{pmatrix}
    ax_1+bx_2\\cx_1+dx_2
  \end{pmatrix}=
  \begin{pmatrix}
    a&b\\c&d
  \end{pmatrix}  \begin{pmatrix}
    x_1\\x_2
  \end{pmatrix}
  , \quad \text{ where } g=\begin{pmatrix}
    a&b\\c&d
  \end{pmatrix}\in\SL.
\end{equation}
The action is linear, that is, it preserves the equivalence relation \(  \begin{pmatrix}
  x_1\\x_2
\end{pmatrix}
\sim   \begin{pmatrix}
  \lambda x_1\\ \lambda x_2
\end{pmatrix}\), \(\lambda\neq 0\) on \(\Space{R}{2}\). The
collection of all respective cosets for non-zero vectors in \(\Space{R}{2}\) is the \emph{projective
  line}%
\index{projective line} \(P\Space{R}{1}\). It is said that a
non-zero vector \(\begin{pmatrix} x_1\\x_2
\end{pmatrix}\in\Space{R}{2}\) represents the point with
homogeneous coordinates%
\index{homogeneous coordinates}
\([x_1:x_2]\in P\Space{R}{1}\). For \(x_2\neq 0\) the point is less tautologically
represented by \([\frac{x_1}{x_2}:1]\) as well. 
The embedding \(\Space{R}{} \rightarrow P\Space{R}{1}\) defined by  \(x \mapsto
[x:1]\), 
\(x\in\Space{R}{}\) covers all but one
point in \(P\Space{R}{1}\). The exception---point \([1:0]\)---is
naturally identified with infinity.

Linear action~\eqref{eq:left-matrix-mult} induces a morphism of
the projective line \(P\Space{R}{1}\),
which is called a \emph{M\"obius transformation}%
\index{M\"obius transformation}.
% This \(\SL\)-action on \(P\Space{R}{1}\) is denoted as \(g: x\mapsto g\cdot x\).
Its restriction to  \(\Space{R}{} \subset P\Space{R}{1}\)  takes linear fractional%
\index{linear fractional transformation} form, cf.~\eqref{eq:linear-fractional-naive}:
\begin{displaymath}
  g:\ 
  [x:1]
  \mapsto
  \left[\frac{ax+b}{cx+d}:1\right]\,, \quad \text{ where } g=\begin{pmatrix}
    a&b\\c&d
  \end{pmatrix}\in\SL \text{ and } cx+d\neq 0.
\end{displaymath}
It is useful to consider some special subgroups of \(\SL\):
\begin{itemize}
\item Compact subgroup \(K\) acts by ``rotations'' without any fixed points
\begin{equation}
  \label{eq:subgroup-K}
  K=\left\{
     \begin{pmatrix}
       \cos t & \sin t \\ -\sin  t& \cos t
     \end{pmatrix}  \such t \in\Space{T}{} 
     \right\}.
\end{equation}
\item Subgroup \(N\) acts by shifts on \(\Space{R}{}\) and has the only fixed point \(\infty=[1:0]\),
\begin{equation}
  \label{eq:subgroup-N}
  N=\left\{
     \begin{pmatrix}
       1 & t \\ 0 & 1
     \end{pmatrix}  \such t \in\Space{R}{} 
     \right\};
\end{equation}
\item Subgroup \(N'\) (conjugated to \(N\)) has the only fixed point \(0=[0:1]\),
\begin{equation}
  \label{eq:subgroup-N-prime}
  N'=\left\{
     \begin{pmatrix}
       1 & 0 \\ t & 1
     \end{pmatrix}  \such t \in\Space{R}{} 
     \right\};
\end{equation}
\item Subgroup \(A\) acts by dilations and fixes both \(0=[0:1]\) and \(\infty=[1:0]\)
\begin{equation}
  \label{eq:subgroup-A}
  A=\left\{
     \begin{pmatrix}
       \rme^t & 0 \\ 0 & \rme ^{-t}
     \end{pmatrix}  \such t \in\Space{R}{} 
     \right\}.
\end{equation}
\end{itemize}

A direct check shows that
\begin{equation}
  \label{eq:preservation-J}
  g^t J g = J, \quad \text{ where } J =
  \begin{pmatrix}
    0 & 1 \\ -1 &0
  \end{pmatrix}, 
\end{equation}
and \(g^t\) is a transpose of \(g  \in \SL\). The above matrix \(J\) defines the \emph{symplectic form}%
\index{symplectic form} \(\omega\) on \(\Space{R}{2}\)~\cite[\S\,41]{Arnold91}:
\begin{displaymath}
  \omega(v_1, v_2) \coloneqq v_1^t J v_2 = x_1 y_2-x_2y_1, \qquad \text{ where } v_i=
  \begin{pmatrix}
    x_i\\y_i
  \end{pmatrix}, \quad i=1,2.
\end{displaymath}
It is a non-degenerate skew-symmetric  bilinear form. In particular \(\omega(v_1,v_2)=0\) if and only if \(v_1\) and \(v_2\) are collinear, it is said that one-dimensional subspaces are \emph{isotropic}\index{isotropic}. Furthermore, non-degeneracy of \(\omega\) implies that one-dimensional subspaces are \emph{maximal} isotropic subspaces, known as \emph{Lagrangian}. Thus, we can say  that points of the projective line \(P\Space{R}{}\)  are presented by Lagrangian subspaces in \(\Space{R}{2}\) and will extend this interpretation to a wider context below.

The identity~\eqref{eq:preservation-J} implies that \(\SL\) preserves \(\omega\), indeed
\begin{equation}
  \label{eq:demo-SL2R=Sp1}
  \omega(gv_1, gv_2)=(gv_1)^t J gv_2 =v_1^t (g^tJ g)v_2 = v_1^t J v_2=\omega(v_1,v_2) 
\end{equation}
for all  \(g \in\SL\). Thus, \(\SL\) coincides with the symplectic group \(\Sp[1]\) in \(\Space{R}{2}\).

\section{Lagrangian subspace of a symplectic space and $\SL$ action}
\label{sec:sympl-struc-lagr}

Let \(\FSpace{H}{}\) be a separable real Hilbert space of finite or infinite dimensionality. Let \(\FSpace{M}{}\) be a closed subspace of \(\FSpace{H}{}\) such that dimensionality of \(\FSpace{M}{}\) and its orthogonal complement \(\FSpace[\perp]{M}{}\)  are equal: \(\dim \FSpace{M}{} = \dim \FSpace[\perp]{M}{}\). Thereafter, we present \(\FSpace{H}{}\) as \(\FSpace{H}{}=\FSpace[\perp]{M}{} \operp \FSpace{M}{}\) (where \(\FSpace{A}{} \operp \FSpace{B}{}\) denotes the direct sum of two \emph{orthogonal} subspaces \(\FSpace{A}{}\) and \(\FSpace{B}{}\)). Correspondingly, we write a generic element \((y,x)\in \FSpace{H}{}\) with \(x\in \FSpace{M}{}\) and \(y\in \FSpace[\perp]{M}{}\). 

Since \(\dim \FSpace{M}{} = \dim \FSpace[\perp]{M}{}\) there exists a unitary operator \(\oper{U}: \FSpace{M}{} \rightarrow \FSpace[\perp]{M}{}\), that is \(\oper{U}\oper{U}^*=\mathbbm{1}_{\FSpace[\perp]{M}{}}\) and \(\oper{U}^*\oper{U}=\mathbbm{1}_{\FSpace{M}{}}\). Let us define a symplectic map \(\oper{W}: \FSpace{H}{}\rightarrow \FSpace{H}{}\) through \(2\times 2\) operator matrix on \(\FSpace{H}{}=\FSpace[\perp]{M}{} \operp \FSpace{M}{}\), cf.~\eqref{eq:preservation-J}:
\begin{equation}
  \label{eq:defn-W-operator}
  \oper{W}=
  \begin{pmatrix}
    \mathbf{0} & \oper{U}\\
    -\oper{U}^* & \mathbf{0}
  \end{pmatrix}.
\end{equation}
Clearly \(\oper{W}^2=-\mathbbm{1}_{\FSpace{H}{}}\), where \(\mathbbm{1}_{\FSpace{H}{}}\) is the identity map on \(\FSpace{H}{}\). Therefore, we define the symplectic form \(\omega: \FSpace{H}{}\times \FSpace{H}{} \rightarrow \Space {R}{}\) by:
\begin{equation}
  \label{eq:defn-symplectic-form}
  \omega((y,x);(y',x')) = \scalar{\oper{W}(y,x)}{ (y',x')}
  =  \scalar[{\FSpace[\perp]{M}{}}]{\oper{U}x}{y'}-\scalar[{\FSpace{M}{}}]{\oper{U}^*y}{x'},
\end{equation}
where \(x, x'\in \FSpace{M}{}\) and \(y,y'\in \FSpace[\perp]{M}{}\). Here we denote by \(\scalar[{\FSpace{M}{}}]{\cdot}{\cdot}\) and \(\scalar[{\FSpace[\perp]{M}{}}]{\cdot}{\cdot}\) the restrictions of the inner product from \(\FSpace{H}{}\) to \(\FSpace{M}{}\) and \(\FSpace[\perp]{M}{}\) respectively. It is easy to show that \(\FSpace{M}{}\) is a \emph{Lagrangian subspace}%
\index{Lagrangian subspace} with respect to \(\omega\), that is \(\FSpace{M}{} = \FSpace[{\sperp}]{M}{}\), where the \emph{symplectic complement}%
\index{symplectic complement} \(\FSpace[{\sperp}]{M}{}\) of \(\FSpace{M}{}\) is:
\begin{displaymath}
  %\label{eq:symplectic-ortho-compl}
  \FSpace[{\sperp}]{M}{} \coloneqq   \{ z \in \FSpace{H}{} \such \omega(z;z')=0 \quad \text{for all} \quad  z' \in \FSpace{M}{}\}.
\end{displaymath}
It is known, that for a finite dimensional \(\FSpace{H}{}\) and any its Lagrangian subspace \(\FSpace{N}{}\) one has \(\dim \FSpace{N}{} = \frac{1}{2} \dim \FSpace{H}{}\). For infinite dimensional situation see the standard reference~\cite{ChernoffMarsden74a} and some recent works~\cites{AbbondandoloMajer15a,Longo22a} with further literature thereafter.

Alternatively, we may start from a Hilbert space \(\FSpace{H}{}\) and a symplectic (i.e., non-degenerate skew-symmetric bilinear) form \(\omega\) on it.  Let \(\FSpace{M}{}\) be a Lagrangian subspace of \(\FSpace{H}{}\), then every non-zero \(z\in \FSpace[\perp]{M}{}\) defines a linear functional \(l_z(\cdot)=\omega(\cdot;z)\) on \(\FSpace{M}{}\), which is  non-zero because \(z\not\in \FSpace[{\sperp}]{M}{}=\FSpace{M}{}\). For similar reason the map \(z \mapsto l_z\) is injective. This linear map is surjective as well: take \(l\in \FSpace[*]{M}{}\) and let \(\FSpace{M}{l}\) be its kernel,  which is a closed proper subspace of \(\FSpace{M}{}\). Therefore,  \(\FSpace[{\sperp}]{M}{} = \FSpace{M}{}\) is a proper subspace of \(\FSpace[{\sperp}]{M}{l}\), which shall contain a non-zero \(z' \in \FSpace[{\sperp}]{M}{l} \setminus \FSpace{M}{}\). Put \(z'' = P_{\FSpace[\perp]{M}{}} z'\), then \(z''\in \FSpace[\perp]{M}{} \cap \FSpace[{\sperp}]{M}{l}\) since \(z''-z' \in \FSpace{M}{} \subset  \FSpace[{\sperp}]{M}{l}\). For some  \(x \in \FSpace{M}{}\) with a non-zero value of \(l(x)\) let us define \(z= \frac{l(x)}{\omega(x;z'')} z'' \in \FSpace[\perp]{M}{}\). Then, \(l_z(\cdot)=\omega(\cdot;z)\) coincides with \(l\) on \(\FSpace{M}{}\). Thus, \(\omega\) defines a bijection \(\FSpace[*]{M}{} \rightarrow \FSpace[\perp]{M}{}\). Combining this with the Riesz--Fr\'echet representation theorem we obtain the following
\begin{proposition}
  Let \(\FSpace{H}{}\) be a Hilbert space equipped with a symplectic form \(\omega\) and  \(\FSpace{M}{}\subset \FSpace{H}{}\) be a Lagrangian subspace. Then, there is a linear bijection
  \begin{equation}
    \label{eq:U-map-M-Mperp-defn}
    \oper{U}:\FSpace{M}{} \rightarrow \FSpace[\perp]{M}{}: x \mapsto y \quad \text{ such that }    \scalar{z}{x} = \omega(z, y) \quad \text{ for all } \quad z \in \FSpace{M}{}.
  \end{equation}
  Therefore, \(\omega\) is represented by~\eqref{eq:defn-symplectic-form} in terms of \(\oper{U}\).
\end{proposition}
In the rest of the paper, we always assume that the above map~\eqref{eq:U-map-M-Mperp-defn} is unitary. Thereafter,  we are avoiding a chicken-egg-type discussion and from the start assume the simultaneous presence of
\begin{itemize}
\item a symplectic form \(\omega\) on \(\FSpace{H}{}\);
\item a Lagrangian subspace \(\FSpace{M}{}\) of \(\FSpace{H}{}\); and
\item a unitary map \(\oper{U}:  \FSpace{M}{} \rightarrow \FSpace[\perp]{M}{}\);
\end{itemize}
which are connected by relations~\eqref{eq:defn-symplectic-form} and~\eqref{eq:U-map-M-Mperp-defn}. In common terms, the choice of Lagrangian \(\FSpace{M}{}\) provides a \emph{polarisation}%
\index{polarisation} of a symplectic space \(\FSpace{H}{}\). Through the rest of this paper we assume a fixed polarisation of \(\omega\) on \(\FSpace{H}{}\)  by a Lagrangian subspace \(\FSpace{M}{}\).
\begin{remark}
  Although it is not directly related to our present consideration, we mention that  symplectic and complex structures are closely connected through the inner product~\cite{Arnold67a}. Furthermore, for a complex inner product space the imaginary part of the inner product is a symplectic form. In the opposite direction: if a real inner product space has a compatible symplectic form, then it can be combined with the real inner product to make a complex inner product, cf.~\cite[\S\,2]{Longo22a}.
\end{remark}

\section{$\SL$-action on a polarised Hilbert space}
\label{sec:sl-action-polarised}

For a given isometric isomorphism \(\oper{U}:\FSpace{M}{}\rightarrow \FSpace[\perp]{M}{}\) and any \(2\times 2 \)-matrix \( g=
\begin{pmatrix} 
a&b\\c&d
\end{pmatrix}\) with scalar entries
 we define a linear operator
\(\oper{M}_{\oper{U}}: \FSpace{H}{}\rightarrow \FSpace{H}{}\)  by:
\begin{equation}
  \label{eq:matrix-action}
  \oper{M}_{\oper{U}}(g)
  =\begin{pmatrix}
     a\mathbbm{1}_{\FSpace[\perp]{M}{}} &b\oper{U}\\
     c\oper{U}^{*}&d\mathbbm{1}_{\FSpace{M}{}}
  \end{pmatrix}:\quad
  \begin{pmatrix}
    y\\x
  \end{pmatrix}
  \mapsto
  \begin{pmatrix}
    ay+b\oper{U}x\\ c \oper{U}^{*}y+dx
  \end{pmatrix} \quad \text{where} \quad
  x\in\FSpace{M}{}, y\in\FSpace[\perp]{M}{}.
\end{equation}
Avoiding matrices the same operator can be written as:
\begin{equation}
  \label{eq:matrix-action-scal}
  \oper{M}_{\oper{U}}(g)=(a\mathbbm{1}_{\FSpace[\perp]{M}{}} +c\oper{U}^{*})P_{\FSpace[\perp]{M}{}} + (d\mathbbm{1}_{\FSpace{M}{}}+b\oper{U}) P_{\FSpace{M}{}},
\end{equation}
where \(P_{\FSpace{M}{}}\) and \(P_{\FSpace[\perp]{M}{}}\) are orthogonal projections of \(\FSpace{H}{}\) on the respective subspaces.
A direct check shows that~\eqref{eq:matrix-action} is a linear
representation of \(\SL\) on \(\FSpace{H}{}\) by bounded operators: \(\oper{M}_{\oper{U}}(g_1)\oper{M}_{\oper{U}}(g_2)=\oper{M}_{\oper{U}}(g_1g_2)\). The
representation is reducible, every two-dimensional subspace spanned by
vectors \(\{x, \oper{U}x\}\) with a non-zero \(x\in \FSpace{M}{}\) is invariant and the
group \(\SL\) acts transitively on every such two-dimensional
plane. In other words, such a plane is an \(\SL\)-homogeneous space. Similarly to~\eqref{eq:demo-SL2R=Sp1} one checks that the symplectic form \(\omega\) is invariant under the action~\eqref{eq:matrix-action}--\eqref{eq:matrix-action-scal}.

\begin{lemma}
\label{lem:homogeneous-space}
  If a non-zero \(z\in \FSpace{H}{}\) is such that
  \begin{enumerate}
  \item    \(\oper{U}(P_{\FSpace{M}{}} z)\) and \(P_{\FSpace[\perp]{M}{}} z\) are collinear, then \(\SL\)-orbit of \(z\) is a two-dimensional subspace  isomorphic to \(\SL/N'\), where \(N'\) is the subgroup~\eqref{eq:subgroup-N-prime} of \(\SL\).
  \item Non-zero vectors \(\oper{U}(P_{\FSpace{M}{}} z)\) and \(P_{\FSpace[\perp]{M}{}} z\) are not collinear, then  \(\SL\)-orbit of \(z\) is a non-linear subset of the four-dimensional vector space spanned by
    \begin{displaymath}
      \{P_{\FSpace{M}{}} z, \oper{U}(P_{\FSpace{M}{}} z), P_{\FSpace[\perp]{M}{}} z, \oper{U}^* (P_{\FSpace[\perp]{M}{}} z)\},
    \end{displaymath}
    which is in bijection with \(\SL\).
  \end{enumerate}
\end{lemma}
\begin{example}
  Let \(\FSpace{H}{}=\Space{R}{4}\) and \((e_n)\), \(n=1,2,3,4\) be the standard orthonormal basis. Let \(\FSpace{M}{}\) be the linear span \(\Lin(e_1,e_2)\) with \(\FSpace[\perp]{M}{}=\Lin(e_3,e_4)\).   The isometry \(\oper{U}: \FSpace{M}{} \rightarrow \FSpace[\perp]{M}{}\) is given by \(\oper{U}(e_1)=e_3\) and \(\oper{U}(e_2)=e_4\). The respective symplectic form is:
  \begin{displaymath}
    \omega( (z_1, z_2, z_3,z_4); (\tilde{z}_1, \tilde{z}_2, \tilde{z}_3,\tilde{z}_4))
    = z_1 \tilde{z}_3 +z_2 \tilde{z}_4 - z_3 \tilde{z}_1 -z_4 \tilde{z}_2.
  \end{displaymath}
  Consider a Lagrangian subspace \(\FSpace{N}{}=\Lin(e_1+e_3, e_2+2e_4)\).
  Therefore, a generic vector in \(\FSpace{N}{}\) has the form \(x(e_1+e_3)+y(e_2+2e_4)\).
  The associated operator \(\oper{T}: \FSpace{M}{} \rightarrow \FSpace{M}{}\) has two eigenvectors \(e_1\) and \(e_2\) with eigenvalues \(1\) and \(2\) respectively. The same eigenvectors deliver so-called \emph{main angles} of the operator sinus of the pair of subspaces \((\FSpace{M}{},\FSpace{N}{})\)~\cite{BoettcherSpitkovsky10a}. 
  
  If we apply \(g=
  \begin{pmatrix}
    a&b\\c&d
  \end{pmatrix} \in \SL\) in \(\Space{R}{4}\) by~\eqref{eq:matrix-action} to \(\FSpace{N}{}\) then we get the following subset of \(\Space{R}{4}\):
  \begin{displaymath}
    % \FSpace[g]{N}{}=
    \{  x(c +  d),  y(c +2  d ), x (a+ b) ,  y(a +2  b )\in \Space{R}{4} \such (x,y) \in \Space{R}{2}  \}.
  \end{displaymath}
  Two observations on the eigenvectors \(e_1\) and \(e_2\) of \(\oper{T}\) and the respective vectors \(e_1+e_3\) and \(e_2+2e_4\) in \(\FSpace{N}{}\):
  \begin{enumerate}
  \item Only scalar multiples of \(e_1+e_3\) or \(e_2+2e_4\) produce a two-dimensional homogeneous space isomorphic to \(\SL/N'\),  where \(N'\) is the subgroup~\eqref{eq:subgroup-N-prime}.
  \item Only scalar multiples of \(e_1+e_3\) or \(e_2+2e_4\) produce orbits which are a homogeneous subsets of \(\Space{R}{4}\).
  \end{enumerate}
%  It will be worth to derive some general theory hinted by observations from this example.
\end{example}

\section{Pairs of Lagrangian subspaces and graphs of operators}
\label{sec:pairs-subsp-comp}

Now we add a second closed subspace \(\FSpace{N}{}\subset \FSpace{H}{}\) and it is common~\cite{BoettcherSpitkovsky10a} to denote the following closed subspaces:
\begin{equation}
  \label{eq:MN-intersections}
  \FSpace{M}{00}=\FSpace{M}{}\cap \FSpace{N}{}, \quad \FSpace{M}{01}=\FSpace{M}{}\cap \FSpace[\perp]{N}{}, \quad \FSpace{M}{10}=\FSpace[\perp]{M}{}\cap \FSpace{N}{}, \quad \FSpace{M}{11}=\FSpace[\perp]{M}{}\cap \FSpace[\perp]{N}{}.
\end{equation}
Also, one introduces closed subspaces \(\FSpace{M}{0} \subset \FSpace{M}{}\) and \(\FSpace{M}{1}\subset \FSpace[\perp]{M}{}\) such that:
\begin{equation}
  \label{eq:M-intersection-decompositions}
  \FSpace{M}{}=\FSpace{M}{00}\operp \FSpace{M}{01} \operp \FSpace{M}{0} \qquad \text{and} \qquad
  \FSpace[\perp]{M}{}=\FSpace{M}{10}\operp \FSpace{M}{11} \operp \FSpace{M}{1} .
\end{equation}
The traditional path~\cites{Halmos69a,BoettcherSpitkovsky10a} is to use to a pair subspaces \(\FSpace{M}{}\) and \(\FSpace{N}{}\) to define a map \(\oper{U}: \FSpace{M}{} \rightarrow \FSpace[\perp]{M}{}\) such that \(\FSpace{N}{}=\{(\oper{U}\oper{T}x, x) \such x \in \FSpace{M}{}\}\) is a graph\footnote{We swap order \((\oper{A}x,x)\) of components in the graph of \(\oper{A}\) in comparison to the traditional one \((x,\oper{A}x)\) to better align with the traditional usage of homogeneous coordinates on the projective line.} of some \emph{self-adjoint} operator \(\oper{T}: \FSpace{M}{} \rightarrow \FSpace{M}{}\). In this way \emph{a pair \((\FSpace{M}{}, \FSpace{N}{})\) of closed subspaces in a generic position defines a symplectic structure on \(\FSpace{H}{}\)}. Since we are approaching from a different end---the symplectic form \(\omega\) and the unitary operator \(\oper{U}\) are already given---we instead put a compatibility condition on  \(\FSpace{N}{}\).
\begin{definition}
  We say that a closed subspace \(\FSpace{N}{}\) of \(\FSpace{H}{}\) is \emph{compatible} with the isometry \(\oper{U}: \FSpace{M}{} \rightarrow \FSpace[\perp]{M}{}\) if \(\FSpace{N}{}\)  is Lagrangian with respect to the symplectic form \(\omega\) from~\eqref{eq:defn-symplectic-form}.
\end{definition}
A compatible subspace \(\FSpace{N}{}\) is specially aligned with the referencing subspace  \(\FSpace{M}{}\): 
\begin{lemma}
  \label{le:subspace-N-map-properties}
  If a closed subspace \(\FSpace{N}{}\) is compatible with the isometry  \(\oper{U}: \FSpace{M}{} \rightarrow \FSpace[\perp]{M}{}\) then
  \begin{enumerate}
  \item \label{item:dim-ker=dim-cocker}
    \(\oper{U} \FSpace{M}{00} = \FSpace{M}{11}\) in particular \(\dim \FSpace{M}{00} = \dim \FSpace{M}{11}\).
  \item \(\oper{U} \FSpace{M}{01} = \FSpace{M}{10}\)  in particular \(\dim \FSpace{M}{01} = \dim \FSpace{M}{10}\).
  \item \(\oper{U} \FSpace{M}{0}=\FSpace{M}{1}\)  in particular \(\dim \FSpace{M}{0}=\dim \FSpace{M}{1}\).
  \end{enumerate}
\end{lemma}
In connection to \(\SL\)-action on operators introduced below, a closed subspace \(\FSpace{N}{}\) can be usefully represented as an orthogonal sum of three closed subspaces:
\begin{displaymath}
  \FSpace{N}{}=\FSpace{N}{\Nnull} \operp \FSpace{N}{\Nodin} \operp \FSpace{N}{\infty},\quad  \text{where} \quad  \FSpace{N}{\Nnull}=\FSpace{M}{00}(=\FSpace{M}{} \cap \FSpace{N}{}), \quad  \FSpace{N}{\infty}=\FSpace{M}{10}(=\FSpace[\perp]{M}{} \cap \FSpace{N}{})
\end{displaymath}
and \(\FSpace{N}{\Nodin}\) is the orthogonal complement of \(\FSpace{N}{\Nnull} \operp \FSpace{N}{\infty}\) within \(\FSpace{N}{}\). Informally, components  \(\FSpace{N}{\Nnull}\) and \(\FSpace{N}{\infty}\) represent eigenspaces for the zero and infinity eigenvalues respectively, while \(\FSpace{N}{\Nodin}\) collects the rest.

We note that \(\FSpace{N}{\Nodin}\subset \FSpace{M}{0} \operp \FSpace{M}{1}\) and two subspaces \(\FSpace{N}{\Nodin}\) and \(\FSpace{M}{0}\) (as well as \(\FSpace{N}{\Nodin}\) and \(\FSpace{M}{1}\)) are in \emph{generic position}~\cite{Halmos69a},  in the sense that all four of the special
intersections of the form~\eqref{eq:MN-intersections} are equal to \(0\) within the Hilbert space  \(\FSpace{M}{0} \operp \FSpace{M}{1}\).
\begin{lemma}
  Let \(\FSpace{N}{}\) be a closed subspace compatible with the isometry  \(\oper{U}: \FSpace{M}{} \rightarrow \FSpace[\perp]{M}{}\).  
  There exits is a (possibly unbounded) closed \emph{self-adjoint} operator \(\oper{T}\) on \(\FSpace{M}{0}\) with the trivial kernel and a dense range such that
  \begin{equation}
    \label{eq:subspace-graph-T}
    \FSpace{N}{\Nodin} = \{(\oper{U}\oper{T}x, x) \such x \in \mathrm{Dom}(\oper{T}) \}. 
  \end{equation}
  The domain and range of \(\oper{T}\) are \(P_{\FSpace{M}{}} \FSpace{N}{\Nodin}\)  and  \(\oper{U}^* P_{\FSpace[\perp]{M}{}} \FSpace{N}{\Nodin}\) respectively.
\end{lemma}
\begin{proof}
  Following the standard arguments (the first paragraph of the proof of Thm.~1 in~\cite{Halmos69a}) one shows that    \(P_{\FSpace{M}{}} \FSpace{N}{\Nodin}\) is dense in \(\FSpace{M}{0}\) and  \(P_{\FSpace[\perp]{M}{}} \FSpace{N}{\Nodin}\) is dense in \(\FSpace{M}{1}\). Then, for \(x\in  P_{\FSpace{M}{}} \FSpace{N}{\Nodin}\) there is a unique point \((y, x) \in \FSpace{N}{}\) and we assign \(\oper{T} x = \oper{U}^* y\). Clearly, \(\oper{U}\oper{T}\) is a linear bijection between \(P_{\FSpace{M}{}} \FSpace{N}{\Nodin}\) and  \(P_{\FSpace[\perp]{M}{}} \FSpace{N}{\Nodin}\) and \(\FSpace{N}{\Nodin}\) is its closed graph as in~\eqref{eq:subspace-graph-T}. By \(\FSpace{N}{} \subset \FSpace[\circ]{N}{} \) of the Lagrangian property of \(\FSpace{N}{}\),  for any \(x_1, x_2\in P_{\FSpace{M}{}} \FSpace{N}{\Nodin}\) and the respective \(y_1=\oper{U}\oper{T} x_1\) , \(y_2=\oper{U}\oper{T} x_2\) we have, cf.~\eqref{eq:defn-symplectic-form}:
  \begin{equation}
    \label{eq:lagrangian-symmetric}
    \begin{split}
     0 &= \omega((x_1,y_1); (x_2,y_2)) \\
       & = \scalar[{\FSpace[\perp]{M}{}}]{\oper{U}x_1}{y_2}-\scalar[{\FSpace{M}{}}]{\oper{U}^*y_1}{x_2}\\
       & = \scalar[{\FSpace[\perp]{M}{}}]{\oper{U}x_1}{\oper{U}\oper{T} x_2}-\scalar[{\FSpace{M}{}}]{\oper{U}^*\oper{U}\oper{T} x_1}{x_2}\\
       & = \scalar[{\FSpace{M}{}}]{x_1}{\oper{T} x_2}-\scalar[{\FSpace{M}{}}]{\oper{T} x_1}{x_2}.
    \end{split}
   \end{equation}
   That is \(\scalar[{\FSpace{M}{}}]{x_1}{\oper{T} x_2}=\scalar[{\FSpace{M}{}}]{\oper{T} x_1}{x_2}\), which shows that \(T\) is a closed symmetric operator. To see that \(\oper{T}\) is self-adjoint we note that if \(\oper{T} \neq \oper{T}^*\) then the graph of \(\oper{T}\) shall be a proper subspace of the closed graph of \(\oper{T}^*\) \cite[\S\,VIII.2]{SimonReed80}. The latter shall be an isotropic subspace (by computation~\eqref{eq:lagrangian-symmetric} read in the opposed direction)---but this contradicts to \(\FSpace[\circ]{N}{} \subset \FSpace{N}{} \) of the Lagrangian property of \(\FSpace{N}{}\).
 \end{proof}
 The above connection~\eqref{eq:subspace-graph-T} between subspaces and graphs of operators can be trivially extended to a  self-adjoint operator \(\hat{\oper{T}}: \FSpace{M}{0}\operp \FSpace{M}{00} \rightarrow  \FSpace{M}{0}\operp \FSpace{M}{00}\) by letting:
 \begin{equation}
   \label{eq:extend-operator-kernel}
   \hat{\oper{T}}(x, x_0)=(\oper{T}x,0), \quad \text{ where } x \in \FSpace{M}{0} \text{ and } x_0 \in \FSpace{M}{00}.
 \end{equation}
 Since \(\dim \FSpace{M}{00} = \dim \FSpace{M}{11}\) by Lem.~\ref{le:subspace-N-map-properties}, a bounded operator \(\hat{\oper{T}}\) is a \emph{Fredholm operator}%
 \index{Fredholm operator} and \(\FSpace{N}{\Nodin}\operp\FSpace{N}{\Nnull}\) corresponds to the graph of \(\oper{U}\hat{\oper{T}}: \FSpace{M}{0}\operp \FSpace{M}{00} \rightarrow  \FSpace{M}{1}\operp \FSpace{M}{11}\). Because r\^oles of \(\FSpace{M}{}\) and \(\FSpace[\perp]{M}{}\) are interchangeable, we have also a self-adjoint operator \(\check{\oper{T}}: \FSpace{M}{1} \operp \FSpace{M}{10} \rightarrow  \FSpace{M}{1} \operp \FSpace{M}{10}\) such that \( \FSpace{N}{\Nodin}\operp \FSpace{N}{\infty}  \)  is the graph of the operator \(\oper{U}^*\check{\oper{T}}:   \FSpace{M}{1} \operp \FSpace{M}{10} \rightarrow  \FSpace{M}{0} \operp \FSpace{M}{01}\). Again, if \(\check{\oper{T}}\) is bounded and \(\dim \FSpace{M}{10}< \infty\) then \(\check{\oper{T}}\)  is a Fredholm operator.

 The above situation is much in line with the standard treatment of the projective line \(P\Space{R}{}\) as a manifold with two charts:
 \begin{displaymath}
   \phi_1: \Space{R}{} \rightarrow P\Space{R}{} : x \mapsto [x:1], \qquad \text{ and}
   \qquad
   \phi_2: \Space{R}{} \rightarrow P\Space{R}{} : x \mapsto [1:x].
 \end{displaymath}
Both charts overlap over \((-\infty,0)\cup (0,\infty)\) and each contains either \(0=[0:1]\) or \(\infty=[1:0]\) as an additional point to their common part. Similarly \(\oper{U}\hat{\oper{T}}\) and \(\oper{U}^*\check{\oper{T}}\) share \(\FSpace{N}{\Nodin}\) as the common portion of their graphs with \(\FSpace{N}{\Nnull}\) and \(\FSpace{N}{\infty}\) being respective supplements in each case.

We take a liberty to encode the above picture by saying that a generic Lagrangian subspace \(\FSpace{N}{}\) represents the graph of a \emph{bi-Fredholm}%
 \index{bi-Fredholm operator} operator \(\FSpace{M}{}\rightarrow \FSpace{M}{}\). Loosely, we can treat \(\FSpace{M}{01}\) as the eigenspace corresponding to the infinite eigenvalue (in the sense of the projective line), cf. ``arithmetic'' with infinity mentioned in Introduction. An accurate description employs a pair of operators \((\hat{\oper{T}},\check{\oper{T}})\) as above to produce two graphs-charts \( \FSpace{N}{\Nodin}\operp \FSpace{N}{\Nnull}  \)  and \( \FSpace{N}{\Nodin}\operp \FSpace{N}{\infty}  \) on \(\FSpace{N}{}\).

\section{$\SL$-action on Lagrangian subspaces and fix subgroups}
\label{sec:sl-action-lagrangian}

For a subspace \(\FSpace{N}{} \subset \FSpace{H}{}\) and any \(g\in\SL\) we denote by \(\FSpace[g]{N}{}\) the image of \(\FSpace{N}{}\):
\begin{equation}
  \label{eq:Ng-image-of-subspace}
  \FSpace[g]{N}{}=\oper{M}_{\oper{U}}(g)\FSpace{N}{}=\{\oper{M}_{\oper{U}}(g)z \such z\in \FSpace{N}{}\}
\end{equation}
where \(\oper{M}_{\oper{U}}(g)\) is given by~\eqref{eq:matrix-action}. Since this \(\SL\)-action preserves the symplectic structure we obtain:
\begin{lemma}
  Let \(\FSpace{N}{}\) be a closed subspace compatible with an isometry \(\oper{U}: \FSpace{M}{}\rightarrow
  \FSpace[\perp]{M}{}\), then \(\FSpace[g]{N}{}\) is a closed subspace compatible with \(\oper{U}\) as well.
\end{lemma}
% Intersections of subspaces~\eqref{eq:MN-intersections} depend on both \(\FSpace{M}{}\) and \(\FSpace{N}{}\), therefore we introducing the notations adjusted for the transformed subspace \(\FSpace[g]{N}{}\):  
% \begin{equation}
%   \label{eq:MNg-intersections}
%   \FSpace[g]{M}{00}=\FSpace{M}{}\cap \FSpace[g]{N}{}, \quad
%   \FSpace[g]{M}{01}=\FSpace{M}{}\cap \FSpace[g\perp]{N}{}, \quad \FSpace[g]{M}{10}=\FSpace[\perp]{M}{}\cap \FSpace[g]{N}{}, \quad \FSpace[g]{M}{11}=\FSpace[\perp]{M}{}\cap \FSpace[g\perp]{N}{}.
% \end{equation}
% Similarly we define \(\FSpace[g]{M}{0}\) and \(\FSpace[g]{M}{1}\) in line with~\eqref{eq:M-intersection-decompositions}. These subspaces for \(\FSpace{N}{e}=\FSpace{N}{}\) with the identity element \(e=
% \begin{pmatrix}
%   1&0\\0&1
% \end{pmatrix}\) will be still denoted as in \eqref{eq:MN-intersections}--\eqref{eq:M-intersection-decompositions}.

Note that if for \(z\in \FSpace{N}{\Nodin}\) the respective vectors \(\oper{U}(P_{\FSpace{M}{}} z)\) and \(P_{\FSpace[\perp]{M}{}} z\) are collinear then \(P_{\FSpace{M}{}} z\) is an eigenvector of the operator \(\oper{T}\) from~\eqref{eq:subspace-graph-T}. Similarly, vectors from \(\FSpace{N}{\Nnull}\) and \(\FSpace{N}{\infty}\) can be treated as eigenvectors corresponding to zero and infinite (\([1:0]\) in homogeneous coordinates) eigenvalues. In view of Lem.~\ref{lem:homogeneous-space} this opens a connection between \(\SL\)-action and spectral properties of \(\oper{T}\). 
% \begin{lemma}
%   Let \(\FSpace{N}{}\) corresponds to an operator \(\tilde{\oper{T}}: \FSpace{M}{0}\operp \FSpace{M}{00} \rightarrow  \FSpace{M}{0}\operp \FSpace{M}{00}\) from Cor.\,\textup{\ref{co:graph-subspace}}. Then, the  subspace \(\FSpace[g]{N}{}\) for \( g=
%   \begin{pmatrix} 
%     a&b\\c&d
%   \end{pmatrix}\) corresponds to a graph of an operator \(\tilde{\oper{T}}_g: \FSpace{M}{} \rightarrow \FSpace{M}{}\) if and only if
%   either
%   \begin{itemize}
%   \item   \(b\neq 0\) and   \(\lambda=-\frac{a}{b}\) is in the resolvent set of \(\tilde{\oper{T}}\); or
%   \item \(b=0\) and \(\dim \FSpace{M}{10}=\dim \FSpace{M}{01} =\dim \FSpace{N}{\infty}=0\).
%   \end{itemize}
% \end{lemma}

% For a bi-Fredholm operator \(\oper{T}: \FSpace{M}{} \rightarrow \FSpace{M}{}\) we will denote by \(\FSpace{N}{\oper{T}}\) the graph of \(\oper{U}\oper{T}\) in \(\FSpace{H}{}=\FSpace{M}{}\oplus \FSpace[\perp]{M}{}\).

\begin{proposition}
  A scalar \(\lambda\in \Space{R}{}\) is an eigenvalue%
  \index{eigenvalue} of a bi-Fredholm operator \(\oper{T}: \FSpace{M}{} \rightarrow \FSpace{M}{}\) if and only if the intersection \( \cap_g \FSpace[g]{N}{}\) is non-zero,  where the intersection is taken for all \(g\) in the subgroup
  \begin{displaymath}
    N_\lambda =\left\{
      \begin{pmatrix}
         1-\lambda t & \lambda ^2 t \\
        -t &  1+\lambda t  
      \end{pmatrix} \such t \in \Space{R}{}
    \right\}.
  \end{displaymath}
  The subgroup \(N\) from~\eqref{eq:subgroup-N} plays the same r\^ole for the infinite eigenvalue. 
\end{proposition}
\begin{proof}
  The necessity of the condition is straightforward. If \(\oper{T}x=\lambda x\), then \(z= ( \lambda \oper{U}x, x)\in \FSpace{N}{}\). A direct check shows that \(\oper{M}_{\oper{U}}(g)z = z\) for all \(g \in N_\lambda\). Thus, \( z \in \cap_g \FSpace[g]{N}{}\) making it non-zero.

  For sufficiency: let \(z \in  \cap_g   \FSpace[g]{N}{}\) for \(g\in N_\lambda \). Since \(z \in \FSpace[g]{N}{}\) all its inverse transformations \(\oper{M}_{\oper{U}}(g^{-1})z\) are in \(\FSpace{N}{}\). % Now we are considering two cases:
  Assume for the moment that \(\lambda=0\). Then, for \(z=(P_{\FSpace[\perp]{M}{}}z, P_{\FSpace{M}{}}z) \in \FSpace[\perp]{M}{} \operp \FSpace{M}{}\), cf.~\eqref{eq:matrix-action}:
  \begin{displaymath}
    \oper{M}_{\oper{U}}(g)z=(P_{\FSpace[\perp]{M}{}}z, -t \oper{U}^*P_{\FSpace[\perp]{M}{}}z +P_{\FSpace{M}{}}z) ), \qquad \text{ for } g\in N_{0}.
  \end{displaymath}
  Therefore, cf.~\eqref{eq:defn-symplectic-form},
  \begin{align*}
    \omega(\oper{M}_{\oper{U}}(g)z ; z)
    &= \scalar[{\FSpace[\perp]{M}{}}]{\oper{U}(-t \oper{U}^*P_{\FSpace[\perp]{M}{}}z +P_{\FSpace{M}{}}z)}{P_{\FSpace[\perp]{M}{}}z}-\scalar[{\FSpace{M}{}}]{\oper{U}^*P_{\FSpace[\perp]{M}{}}z}{P_{\FSpace{M}{}}z} \\
    &= \scalar[{\FSpace[\perp]{M}{}}]{-t P_{\FSpace[\perp]{M}{}}z +\oper{U}P_{\FSpace{M}{}}z}{P_{\FSpace[\perp]{M}{}}z}-\scalar[{\FSpace{M}{}}]{\oper{U}^*P_{\FSpace[\perp]{M}{}}z}{P_{\FSpace{M}{}}z} \\
    &= -t \scalar[{\FSpace[\perp]{M}{}}]{P_{\FSpace[\perp]{M}{}}z}{P_{\FSpace[\perp]{M}{}}z} +\scalar[{\FSpace[\perp]{M}{}}]{\oper{U}P_{\FSpace{M}{}}z}{P_{\FSpace[\perp]{M}{}}z}-\scalar[{\FSpace{M}{}}]{\oper{U}^*P_{\FSpace[\perp]{M}{}}z}{P_{\FSpace{M}{}}z} \\
    &= -t \norm{P_{\FSpace[\perp]{M}{}}z}^2 .
  \end{align*}
  Since the symplectic form of two vectors from the Lagrangian subspace \(\FSpace{N}{}\) is zero we get \(P_{\FSpace[\perp]{M}{}}z=0\), which means \(z\) is the eigenvector with the eigenvalue \(\lambda=0\).

  For a generic \(\lambda\in\Space{R}{}\) we note that
  \begin{displaymath}
      \begin{pmatrix}
         1-\lambda t & \lambda ^2 t \\
        -t &  1+\lambda t  
      \end{pmatrix}
      =
      \begin{pmatrix}
         c & s \\
        -s &  c  
      \end{pmatrix}
      \begin{pmatrix}
         1 & 0 \\
        -t/cs &  1  
      \end{pmatrix}
      \begin{pmatrix}
         c & -s \\
        s &  c  
      \end{pmatrix},
    \end{displaymath}
    where \(\lambda=\frac{s}{c}\) and \(c^2+s^2=1\). That is, \(N_\lambda\) is conjugated to \(N'\) from~\eqref{eq:subgroup-N-prime} by the element \(      g_{\lambda}= \begin{pmatrix}
         c & -s \\
        s &  c  
      \end{pmatrix}\) from the compact subgroup \(K\) from~\eqref{eq:subgroup-K}, which ``rotates'' \(\lambda\) to \(0\). Also we can directly check that \(P_{\FSpace[\perp]{M}{}}z = \lambda \oper{U}P_{\FSpace{M}{}}z\) if and only if \( P_{\FSpace[\perp]{M}{}} \oper{M}_{\oper{U}}(g_{\lambda})z =0\). Thus, the conjugation with \(\oper{M}_{\oper{U}}(g_{\lambda})\) convert the case of generic \(\lambda\) to the zero eigenvalue considered above.
\end{proof}

% A generic subgroups which fixes the two given points \(\lambda\) and  \(\mu\) is given by 
% \begin{displaymath}
%   \begin{pmatrix}
%     -\mu c & \lambda a \\ -c & a 
%   \end{pmatrix}
%   \begin{pmatrix}
%     \tau & 0 \\  0 & \tau^{-1}
%   \end{pmatrix}
%   \begin{pmatrix}
%     a & -\lambda a \\ c & -\mu c
%   \end{pmatrix}
%   =
%   \frac{1}{{(\lambda-\mu)} \tau}
%       \begin{pmatrix}
%         {\lambda- \tau^{2} \mu}&
%                                     {\lambda \mu(  \tau^{2}-1)} \\
%         {1-\tau^{2}}&
%                          { \lambda \tau^{2}-\mu}
%       \end{pmatrix}
% \end{displaymath}
% where we use the substitution \(a= {(\lambda c + d)^{-1}}\) in the right-hand side to  get the unit determinant.

For some \(\varepsilon \in(0,1)\) and a subspace \(\FSpace{N}{} \subset \FSpace[\perp]{M}{} \operp \FSpace{M}{}\)  we define its \(\varepsilon\)-adjunct neighbourhood \(\FSpace{N}{\varepsilon}\) by
\begin{equation}
  \label{eq:adjunct-neighborhood}
  \FSpace[\varepsilon]{N}{} =\{ \, z+y \such z\in \FSpace{N}{} , \ y\in \FSpace[\perp]{M}{} \text{ and }  \norm{y} < \varepsilon \norm{P_{\FSpace{M}{}}z} \, \} \, .
\end{equation}
In the model two dimensional case of \(\FSpace[\perp]{M}{} \operp \FSpace{M}{}=\Space{R}{2}\) and \(\FSpace{N}{}\) being a line \(y=\lambda x\) for some \(\lambda\in \Space{R}{}\)   its \(\varepsilon\)-adjunct neighbourhood \(\FSpace[\varepsilon]{N}{}\) consists of all points \((\mu x, x)\) with \(\modulus {\mu-\lambda} < \varepsilon\).

\begin{proposition}
  Let \(\FSpace{N}{}\) be a Lagrangian subspace with the trivial \(\FSpace{N}{\infty}=\FSpace{M}{}\cap \FSpace[\perp]{N}{}\). Let \(\oper{T}: \FSpace{M}{} \rightarrow \FSpace{M}{}\) be the corresponding operator with the graph represented by \(\FSpace{N}{}\). The operator \(\oper{T}-\lambda\mathbbm{1}\) for \(\lambda\in\Space{R}{}\) does not have a bounded inverse provided for a fixed \(\varepsilon \in (0,1)\)  the intersection \(\FSpace[\varepsilon]{N}{} \cap \FSpace[g]{N}{}\) is non-zero for any \(g\) from the semi-group
  \begin{equation}
    \label{eq:sungroup-A-lambda}
    % \frac{ c a \lambda}{t}+ t d a&\frac{ d a \lambda}{t}- t d a \lambda\\\frac{ c a}{t}- t c a& t c a \lambda+\frac{ d a}{t}
    \mathsf{A}_\lambda = \left\{
      \begin{pmatrix}
        \rme^{t}& \lambda (\rme^{-t} - \rme^{t} )\\0&\rme^{-t} 
      \end{pmatrix} \such t \in (0,\infty)
    \right\}.
  \end{equation}
\end{proposition}
\begin{proof}
  For given \(\lambda\in \Space{R}{}\), the subgroup which fixes both \(\lambda\) and  \(\infty\) consists of matrices
  \begin{displaymath}
    \begin{pmatrix}
      \rme^{t} & \lambda (\rme^{-t} - \rme^{t} )\\0&\rme^{-t}
    \end{pmatrix}
    =
    \begin{pmatrix}
      1 & \lambda  \\ 0 & 1 
    \end{pmatrix}
    \begin{pmatrix}
      \rme^{t} & 0 \\  0 & \rme^{-t}
    \end{pmatrix}
    \begin{pmatrix}
      1 & -\lambda  \\ 0 & 1
    \end{pmatrix}
    .
  \end{displaymath}
It is the subgroup \(A\) from~\eqref{eq:subgroup-A} conjugated by the shift \(\lambda\mapsto 0\) from the subgroup \(N\) in~\eqref{eq:subgroup-N}. This time we are shifting the eigenvalue by the element of \(N\) (rather than rotating it by an element of \(K\)) because we want to fix the infinity---the common fixed point of \(A\) and \(N\). As in the proof of the previous Proposition this allows us to work with technically simpler case of \(\lambda=0\) and we will do it from now. The corresponding semi-group~\eqref{eq:sungroup-A-lambda} is identified as \(A_0\) and it is a subset of subgroup \(A\) from~\eqref{eq:subgroup-A}.

A simple observation: if for \(z_1\), \(z_2\in\FSpace{N}{\Nodin}\operp \FSpace{N}{\Nnull}\) we have \(P_{\FSpace{M}{}} z_1 = P_{\FSpace{M}{}} z_2\) then \(z_1=z_2\). For \(g_t\in A_0\) and \(z=(y,x)\in \FSpace{N}{}\) we have \(z_t\coloneqq \oper{M}_{\oper{U}}(g_t)z= (\rme^{t}y, \rme^{-t}x) \in \FSpace[{g_t}]{N}{}\). Since \(\rme^{-t}z=(\rme^{-t}y, \rme^{-t}x)\) is in \(\FSpace{N}{}\), the requirement  of \(z_t\) being also in \(\FSpace[\varepsilon]{N}{}\) requires
\begin{displaymath}
  \norm{\rme^{t}y-\rme^{-t}y}\leq \varepsilon \norm{\rme^{-t}x}, \quad\text{ therefore } \quad
   (\rme^{2t}-1) \norm{y} \leq \varepsilon \norm{x}.
\end{displaymath}

Since  for any \(t\in (0,\infty)\) there exists \((y,x) \in  \FSpace{N}{}\)  such that \(\norm{y} \leq \frac{\varepsilon }{\rme^{2t}-1} \norm{x}\) we see that an inverse of \(\oper{T}\), which shall map \(\oper{U}^*y \) to \(x\),  cannot be bounded.
\end{proof}
In terms of standard spectral theory, the above proposition describes a point \(\lambda\) of the \emph{continuous spectrum}%
\index{continuous spectrum}. 

\section{Conclusion}
\label{sec:conclusion}

Using Halmos' hint from the epigraph, we introduce the following framework:
\begin{itemize}
\item   For a given Hilbert space \(\FSpace{H}{}\) equipped with a symplectic form the associated \emph{operator projective line}%
  \index{operator projective line}%
  \index{projective line|operator} \(P\FSpace{H}{}\) is the family of all Lagrangian subspaces of \(\FSpace{H}{}\). In the finite dimension case, it is known as \emph{Lagrangian Grassmannian}%
  \index{Lagrangian Grassmannian}~\cite{Arnold67a}. 
\item Every Lagrangian subspace \(\FSpace{M}{}\subset \FSpace{H}{}\) defines \(\SL\)\emph{-action} \(\oper{M}_{\oper{U}}(g)\) from~\eqref{eq:matrix-action} on  \(P\FSpace{H}{}\).
\item Let a bi-Fredholm operator \(\oper{T}\) have a Lagrangian subspace \(\FSpace{N}{} \in P\FSpace{H}{}\) as its graph. \emph{Linear fractional transformation}%
\index{operator linear fractional transformation}%
\index{linear fractional transformation!operator}~\eqref{eq:linear-fractional-naive} of \(\oper{T}\) for \(g\in\SL\) is represented by the action \(\oper{M}_{\oper{U}}(g)\) on \(\FSpace{N}{}\).
\end{itemize}
The construction is done within an extension of the \emph{Erlangen programme}%
\index{Erlangen programme} to \emph{non-commutative geometry}%
\index{non-commutative geometry}~\cite{Kisil97a,Kisil02c}.
This setting seems to deserve some further exploration. For example, it will be interesting to the consider operator version of the \emph{Poincare extension}%
\index{Poincare extension} from~\cite{Kisil15a}, which is constructed by a pair of the Lagrangian subspaces \(\FSpace{N}{1}\) and  \(\FSpace{N}{2}\) and their \(\SL\)-transformations defined by a Lagrangian subspace \(\FSpace{M}{}\). A more ambitious aim is to develop a non-commutative (operator) version of the extended M\"obius--Lie \emph{sphere geometry}%
\index{sphere geometry}~\cite{Kisil19a}. There are also numerous inspiring hints on connections to conformal nets%
\index{conformal nets}~\cites{CarpiHillierLongo15a,Longo22a},  algebraic quantum field theory%
\index{algebraic quantum field theory}, and  quantization in general~\cite{Arnold67a}.

\section*{Acknowledgments}
\label{sec:acknowledgments}
The second named author is grateful to Prof.\,I.M.\,Spitkovsky for his long-standing crucial support. The anonymous referee made numerous helpful comments and suggestions which were used to correct and improve the paper.

  \providecommand{\noopsort}[1]{} \providecommand{\printfirst}[2]{#1}
  \providecommand{\singleletter}[1]{#1} \providecommand{\switchargs}[2]{#2#1}
  \providecommand{\irm}{\textup{I}} \providecommand{\iirm}{\textup{II}}
  \providecommand{\vrm}{\textup{V}} \providecommand{\cprime}{'}
  \providecommand{\eprint}[2]{\texttt{#2}}
  \providecommand{\myeprint}[2]{\texttt{#2}}
  \providecommand{\arXiv}[1]{\myeprint{http://arXiv.org/abs/#1}{arXiv:#1}}
  \providecommand{\doi}[1]{\href{http://dx.doi.org/#1}{doi:
  #1}}\providecommand{\CPP}{\texttt{C++}}
  \providecommand{\NoWEB}{\texttt{noweb}}
  \providecommand{\MetaPost}{\texttt{Meta}\-\texttt{Post}}
  \providecommand{\GiNaC}{\textsf{GiNaC}}
  \providecommand{\pyGiNaC}{\textsf{pyGiNaC}}
  \providecommand{\Asymptote}{\texttt{Asymptote}}
\providecommand{\url}[1]{{#1}}
\providecommand{\urlprefix}{URL }
\expandafter\ifx\csname urlstyle\endcsname\relax
  \providecommand{\doi}[1]{DOI~\discretionary{}{}{}#1}\else
  \providecommand{\doi}{DOI~\discretionary{}{}{}\begingroup
  \urlstyle{rm}\Url}\fi

\ifspringer
%\bibliographystyle{spmpsci}
%\bibliography{abbrevmr,akisil,analyse,algebra,arare,aclifford,aphysics,acompute,ageometry,acombin}

\else
  \small %\bibliography{abbrevmr,akisil,analyse,algebra,arare,aclifford,aphysics,acompute,ageometry,acombin}

\fi

\end{document}
\newpage

Depending on the choice of scalars, \(2\times 2\) matrices with
non-zero determinant form groups \(\mathrm{GL}_2(\Space{R}{})\) or
\(\mathrm{GL}_2(\Space{C}{})\), where the operation is given by matrix
multiplication. By the nature of matrix
formula~\eqref{eq:matrix-action-E} the correspondence \(g\mapsto
\oper{M}_g\) is a representation of these groups by linear operators
on \(\FSpace{H}{}\).

A direct check shows that: 
\begin{displaymath}
  \begin{pmatrix}
    I & \oper{T}\\\oper{T}& \oper{T}^2
  \end{pmatrix}=
  \frac{1}{2}
  \begin{pmatrix}
    I & I\\\oper{T}& \oper{T}
  \end{pmatrix}
  \begin{pmatrix}
    I & \oper{T}\\I& \oper{T}
  \end{pmatrix}=
  \frac{1}{2}
  {\begin{pmatrix}
    I & \oper{T}\\I& \oper{T}
  \end{pmatrix}\!}^*
  \begin{pmatrix}
    I & \oper{T}\\I& \oper{T}
  \end{pmatrix}.
\end{displaymath}
Thus, we obtain:
\begin{displaymath}
  {\begin{pmatrix}
    d &-b\\-c&a
  \end{pmatrix}\!}^*
  \begin{pmatrix}
    I & \oper{T}\\\oper{T}& \oper{T}^2
  \end{pmatrix}
  \begin{pmatrix}
    d &-b\\-c&a
  \end{pmatrix}
  \sim
  \begin{pmatrix}
    I & \oper{T}_1\\\oper{T}_1& \oper{T}_1^2
  \end{pmatrix},
\end{displaymath}
where
\begin{displaymath}
  \oper{T}_1=\frac{a\oper{T}-bI}{-c\oper{T}+dI}.
\end{displaymath}

Two subspaces define endpoints of an interval \([x,y]\) on the
projective line. Then an equivalence relation \([x,y]\sim [x',y']\)
may be defined by value of one of the \(\SL\)-invariant bilinear
forms:
\begin{displaymath}
  \begin{pmatrix}
    a&-b\\b&a
  \end{pmatrix}
  \quad
  \begin{pmatrix}
    a&-b\\b&-a
  \end{pmatrix}
  \quad
  \begin{pmatrix}
    a&-b\\b&0
  \end{pmatrix}
\end{displaymath}

\end{document}